\documentclass[11pt,twoside]{amsart}
\usepackage{latexsym,amsfonts,amsmath,amssymb,soul}

\numberwithin{equation}{section}

\makeatletter
\@namedef{subjclassname@2020}{%
  \textup{2020} Mathematics Subject Classification}
\makeatother

\catcode`\@=11

\renewcommand{\section}%
   {\setcounter{equation}{0}\@startsection {section}{1}{\z@}{-3.5ex plus -1ex
  minus -.2ex}{2.3ex plus .2ex}{\Large\bf}}

\newtheorem{Th}{Theorem}[section]
\newtheorem{Rem}[Th]{Remark}
\newtheorem{Rems}[Th]{Remarks}

\newtheorem{Lemma}[Th]{Lemma}

\usepackage{calrsfs}
\DeclareMathAlphabet{\pazocal}{OMS}{zplm}{m}{n}

\def\R{\mathbb R}

\def\sp{\hskip 0.5 pt}
\def\spp{\hskip 1pt}

\def\spi{\hskip 2pt}

\newcommand{\beqsn}{\arraycolsep1.5pt\begin{eqnarray*}}
\newcommand{\eeqsn}{\end{eqnarray*}\arraycolsep5pt}
\newcommand{\beqs}{\arraycolsep1.5pt\begin{eqnarray}}
\newcommand{\eeqs}{\end{eqnarray}\arraycolsep5pt}

\title{A third-order conservation law for the Kirchhoff-Pokhozhaev equation}

\author[Boiti]{Chiara Boiti}
\address{
Dipartimento di Matematica e Informatica \\Universit\`a di Ferrara\\
Via Ma\-chia\-vel\-li n.~30\\
I-44121 Ferrara\\
Italy}
\email{chiara.boiti@unife.it}

\author[Manfrin]{Renato Manfrin}
\address{
Dipartimento di Culture del Progetto, Universit\`a IUAV di Venezia\\
Dorsoduro n.~2196, I-30123 Venezia, 
 Italy}

\email{manfrin@iuav.it}

\begin{document}

\subjclass[2020]{35L65, 35L72, 35L15, 35L20}
\keywords{Kirchhoff equation, third order conservation law.}

\begin{abstract}
We prove that the special Kirchhoff equation studied by Pokhozhaev in \cite{P1, P3} admits  
a third-order conservation law. We further show that if the energy of the solution  is sufficiently small,
then the $L^2$-norms of the derivatives up to third order of the solution remain uniformly bounded with respect to time.
 \end{abstract}

\maketitle

\markboth{\sc  A third-order conservation law for the Kirchhoff-Pokhozhaev equation}
 {\sc C.~Boiti, R.~Manfrin}


\section{Introduction and results}

Let $\Omega \subset \R^n_x$ be a bounded, smooth domain or $\Omega = \R^n_x\sp$. Given a continuous function 
$m:[0,\infty) \rightarrow (0,\infty)$, let us consider the quasi-linear Kirchhoff wave equation 
  \begin{equation} \label{GK}
u_{tt} - m\left (\int_\Omega |\nabla u |^2 dx \right ) \Delta u=0 \quad \text{in} \quad \Omega \times [0,T),
\end{equation}
 for some $T> 0$, with the boundary condition
\begin{equation}\label{bc}
u=0 \quad \text{on} \quad  \partial \Omega \times [0,T),
\end{equation}
if $\sp \Omega\sp  $ is a bounded, smooth domain. As it is known, equation \eqref{GK} is of variational type and each  solution has a  conserved energy. Namely, if $\spp u\spp $ is a sufficiently regular solution of 
\eqref{GK}-\eqref{bc}, then the first-order functional
\begin{equation}\label{CL}
I_1 :=\int_\Omega |u_t |^2 dx + M \left  ( \int_{\Omega} |\nabla u|^2 dx \right ) ,
\end{equation}
where $M(h)= \int_0^h m(s) ds\sp$, remains constant in $\sp [0,T)\sp$. For the Kirchhoff equation, in the general form \eqref{GK}, this is 
essentially the only known conservation law. 

In 1974 Pokhozhaev 
  presented in \cite{P1} a second-order conservation low for the   special
  Kirhhoff equation
   \begin{equation}\label{KP}
  u_{tt}- \frac {\Delta u} {\big ( a\int_{\Omega} |\nabla u |^2 dx + b \big )^2 }  = 0 
\quad  \text{in} \quad \Omega \times [0,T),
   \end{equation}
with  $a, b \in \R \sp $ not both zero. More precisely: 

{\it Let  $\sp u \sp $ be a sufficiently regular function in $\sp \Omega \times [0,T)\sp $ such that
\begin{equation}
a\int_{\Omega} |\nabla u |^2 dx  +b \ne 0 \quad \text{for} \quad t \in [0,T);
\end{equation}
if $\sp u\sp $ is a solution of  \eqref{KP}, \eqref{bc}  then the second-order functional
\begin{equation}\label{CL2}
I_2 \, := \,  ( a \| \nabla u \|^2 + b) \| \nabla u_t\|^2 +
 \frac { \|\Delta u \|^2} 
{ a \|\nabla u \|^2  + b} 
 - {a} \left ( \int_{\Omega} \nabla u \cdot \!\nabla u_t\spp dx \right )^2,
\end{equation}
where $\, \|\cdot \|^2=  \int_{\Omega} |\cdot |^2dx$, \spp remains constant in $[0,T)$.}

Still in  \cite{P1}, Pokhozhaev  proved the remarkable fact that \eqref{KP} is  the unique case in which the Kirchhoff 
equation \eqref{GK} admits 
a second order conservation law. That is, \eqref{GK}
 has a second-order conservation law if and only if 
\begin{equation}
m(s) = \frac{1}  {(a  s +b )^2}\sp ,
\end{equation}
for some $a,b\in \R$ not both zero. See also \cite{BM} for a simple derivation of the expression \eqref{CL2}.

The existence of a conservation law for a second-order functional such as  \eqref{CL2} is important because it implies, under suitable assumptions, a-priori estimates for $\frac{d}{dt}\|\nabla u\|^2$.
In turn, these a-priori estimates lead, as shown in \cite{AP}, 
to the proof of the global solvability of the initial boundary value problems and Cauchy problems for \eqref{KP}. For example, a-priori estimates for $\frac{d}{dt}\|\nabla u\|^2$ can be easily obtained from \eqref{CL2} in the following cases:

\medskip
$(i)$ If $a, b>0$, \eqref{CL2} implies
\begin{equation}
\label{18}
\|\nabla u_t\|^2  \le I_2 b^{-1}
\end{equation}
since $\spp a\|\nabla u\|^2\cdot\|\nabla u_t\|^2 - a (\sp \int_\Omega\nabla u\cdot \!\nabla u_t \spp dx )^2 \ge 0 \spp$. Then, for $t\geq0$,
 $$\|\nabla u(t)\|\le  \|\nabla u(0)\|+  \sqrt { {I_2} b^{-1}}\spp t.$$
  Therefore 
$$\, \left|\frac d {dt}\, \| \nabla u \|^2\right| \le  2 \sqrt {I_2b^{-1}} 
\left  (\|\nabla u(0)\|+  \sqrt { {I_2} b^{-1}}\spp t \right ). $$
This inequality was used in \cite{P3} and then in  \cite{PE} to prove global solvability.

\medskip
$(ii)$ If $a<0$, $b>0$  and $\, a\|\nabla u\|^2+ b > 0
\spp$ (i.e. $ \|\nabla u \|^2 < b |a|^{-1}$), then by  \eqref{CL2} 
\beqsn
I_2\geq|a|\left(\int_\Omega\nabla u\cdot\nabla u_t\right)^2=
|a|\left(\frac12\frac{d}{dt}\|\nabla u\|^2\right)^2.
\eeqsn
Therefore, we have
\begin{equation}
\left |\frac d {dt} \|\nabla u\|^2 \right | \le  2 \sqrt{{I_2} {|a|^{-1}}}.
\end{equation}
This case was considered in \cite {PN}.

In both cases $(i)$ and $(ii)$, the boundedness of $\frac d {dt} \|\nabla u\|^2 $, together with the conservation law for $\sp I_1\sp $, allows the extension of the local solution beyond any fixed $T>0$, leading to a global solution. Moreover, in case $(ii)$,  the global
existence of solutions was established under weak regularity assumptions on the initial data,
lower than those required by the functional $I_2$ (see \cite{PN}).

\medskip
In this paper, we  show that \eqref{KP} has also a third-order conservation law (see
Theorem~\ref{T1}) and then we prove, under the assumption of small initial data, the boundeness of the $L^2$\spp -\spp norms of the  derivatives of the solutions up to third order (see Theorem~\ref{Teo4}).

More precisely, choosing $\Omega=\R^n$ for simplicity, setting
\begin{equation} \label{qs}
s= s(t) := \| \nabla u  \|^2=\int_{\R^n} | \nabla u (x,t)|^2 dx ,
\end{equation}
we prove the following third-order conservation law for the Kirchhoff-Pokhozhaev equation
  \begin{equation}\label{KPR}
  \displaystyle{u_{tt}- \frac {\Delta u} {\big ( a\int_{\R^n} |\nabla u |^2 dx + b \big )^2 }  = 0} \quad 
  \; \text{in} \quad \R^n \times [0,T):
\end{equation}
\begin{Th}\label{T1}
Let $u$ be a function  defined in $\R^n \times [0,T)$, for some $T> 0$, and  suppose 
\begin{equation}\label{reg}
u \in C^k([0,T); H^{3-k}(\R^n)), \quad \text{for} \quad k= 0, \spp 1 
\end{equation}
and 
\begin{equation}\label{reg2}
q = q(t)  = a\sp \|\nabla u \|^2  +b \ne 0 \quad \text{in} \quad [0,T).
\end{equation}
If $\sp u $ is  a solution of \eqref{KPR}, then the third-order functional
\begin{equation}
\label{I3}
\begin{aligned}
I_3\spp &:= \spp q \spp {\|\Delta u_t\|^2} + \frac{\|\nabla (\Delta u)\|^2} q - {q'} \!\int_{\R^n} \Delta u \sp\Delta u_t\spp dx 
 \\
 &\quad\quad \; + \frac{1} 8\spp  {q'}^2 \spp 
  \left ( { q \spp \| \nabla u_t\|^2}  +  \frac{\| \Delta u\|^2} {q}\right )
  -\frac{ a }{16}  \left (  \frac {a^2\sp {s'}^4} 4+ {q^2\spi {s''}^2} \right ) 
 \end{aligned}
 \end{equation}
 remains constant in $[0,T)$.
 \end{Th}

 \begin{Rem} In \eqref{I3} $s'$ and $s''$  are the derivatives of $\spp s(t)$  defined in \eqref{qs}. Under the assumptions
 of Theorem~\ref{T1}, we easily have
 \beqs
 \label{s's''}
 s'=2\int_{\R^n}\nabla u\cdot \!\nabla u_t\sp\sp dx \;  \quad\text{and} \;\quad s'' = 2 \|\nabla u_t\|^2 - \frac 2{q^2}\| \Delta u\|^2,
 \eeqs
where the expression of $\sp s''$ is obtained taking into account that $\spp u$ is a solution of \eqref{KPR}.
 \end{Rem}

Note that, while the first-order functional $I_1$ is always  positive definite (i.e. $ > 0$ if 
$\sp u\sp $ is not constant), 
 $I_2$ and $I_3$ are positive definite
only under some additional assumptions, such as $(i)$ or $(ii)$ above. With the following two lemmas, we prove that the functionals $I_2$ and $I_3$ are positive definite  when $\sp a \ne 0\sp$, $\sp b > 0\sp$
and  $ I_1 $ is sufficiently small. 

More precisely, let us consider  $\sp I_1\sp$ in the case of the 
Kirchhoff-Pokhozhaev equation \eqref{KPR} with  $\sp a \ne 0\sp$ and $\sp b > 0 \sp $. 
That is, according to \eqref{CL}, the functional
\begin{equation}\label{I1KP}
I_1\, := \,  \| u_t\|^2 + \int_0^{\|\nabla u \|^2}  \frac {ds} {(a\sp s + b)^2} \,  =  \,\| u_t\|^2 + 
\frac 1 b \spp \frac {\|\nabla u \|^2} {a\spp \|\nabla u \|^2 + b} \spp ,
\end{equation}
where, obviously, we must assume 
\begin{equation}
\|\nabla u \|^2 < \frac b {|a|} \quad \text{if} \quad a< 0 \sp .
\end{equation} 

\begin{Lemma} 
\label{Teo3} 
Let $\sp a \ne 0 \sp$, $\sp  b  >0 \sp$. Then, we have that
\begin{equation}\label{condd}
I_1 \le \frac 1 {6 |a\sp b|} \;\quad
\Rightarrow \; \quad \frac {b} 2 \le q \le \frac32 \spp b\quad \text{and} \quad |a| \lambda \le \frac 1 2,
\end{equation}
for all $\sp t \in [0,T)$, where
\begin{equation}\label{lambd}
\lambda \, := \, q \spp \| u_t\|^2  + \frac{ \| \nabla u\|^2 }q \spp .
\end{equation}
\end{Lemma}
Furthermore, when $\sp \lambda\sp$ is small and $\sp  q> 0 \sp $, we also have:
\begin{Lemma}\label{Teo2}
Let the assumptions of Theorem~\ref{T1} hold, and also assume that  $q> 0$ in $[0,T)$. 

If $\sp u $ is  a solution of \eqref{KPR}, then
the functionals $\sp I_2, I_3 \sp $ satisfy the inequalities:
\begin{equation}\label{stima0}
\frac 23 \, I_2  \,\le \,  q \spp {\|\nabla u_t\|^2} + \frac{\|\Delta u\|^2} q \, \le \,  2\, I_2, 
\quad\text{if} \quad |a |\lambda \le 2\sp ;
\end{equation}
\begin{equation}
 \label{stima}
\frac 2 3 \,  I_3  \, \le \,   q \spp {\|\Delta u_t\|^2} + \frac{\|\nabla (\Delta u)\|^2} q \, \le \,  2\, I_3 ,
\quad\text{if} \quad |a |\lambda \le \frac 1 2\spp .
\end{equation}
\end{Lemma}

The previous two lemmas  now lead  to the following result:
\begin{Th} \label{Teo4}
Let  $\sp a \ne 0 \sp$, $\sp  b  >0 \sp$. 
Let $u$ be a solution of \eqref{KPR} that satisfies \eqref{reg} and such that
\begin{equation}
I_1  \le \frac{ 1} {6 |a\sp b |}\spp  .
\end{equation} 
Then  $\spp \frac {b}2 \le q \le \frac32 \sp b\spp$  and the estimates \eqref{stima0}, \eqref{stima} are verified 
 for all $\spp t\in [0,T)$. 

In particular, this means that the $L^2$-norms 
$\spp \|\partial_x^{\alpha} \partial_t^ku\|\spp $, for $\spp 1\le |\alpha|+ k \le 3\spp $,
 are uniformly bounded in $[0,T)$, independently of $\spp T> 0$.
\end{Th}

\smallskip
\noindent
\begin{Rems}

In the last statement of Theorem~\ref{Teo4}, concerning  the boundedness of the $L^2$-norms  of  the derivatives of the solution, we can simply assume  $\sp a b \ne 0\sp $. In fact, the equation \eqref{KPR} does not change if we replace $\sp a \sp $ with $\sp -a\sp$  and $\sp b \sp$ with $\sp -b\sp$, respectively. We have used the condition $\sp a\ne 0\sp$, $\sp b > 0\sp$ only for simplicity in the statements and proofs of Lemma~\ref{Teo3} and Theorem~\ref{Teo4}. 

\smallskip
Furthermore, in Theorem~\ref{Teo4}, we have not considered the cases $\sp a=0\sp$, $\sp b\ne 0 \sp$  and $\sp a\ne 0 \sp $, $\sp b= 0 \sp $. 
The first is trivial, since \eqref{KPR} becomes  a linear wave equation. 
The second is, in a certain sense, a degenerate case because when $\sp b=0\sp $ equation \eqref{KPR}  does not have the solution $\sp u\equiv 0 \sp $.

\smallskip

Note that equation \eqref{KPR} keep the regularity of the data also on the solution, but Theorem~\ref{Teo4} gives extra information about boundedness of the derivatives
up to order 3. This boundedness property is not implied by the conservation of the functionals $\sp I_1\sp$ and $\sp I_2\sp$.

\smallskip
Finally, we note that all the previuos results remain valid if $\,\Omega\subset\R^n_x$
is a bounded smooth domain, replacing the Fourier transform with the
Fourier series expansion in eigenfunctions of the Laplace operator (as in \cite{P3}) and adopting the appropriate 
boundary conditions.

\end{Rems}

\medskip
The paper is organized as follows: 
we apply the partial Fourier transform with respect to the $x$-variable in \eqref{KPR} and obtain the Liouville-type second order ordinary differential equation \eqref{lin2}.
In Section~\ref{LL}, we construct a quadratic form for equation \eqref{lin2} that will be applied, in Section~\ref{sec3}, to obtain a conservation law for equation \eqref{KPR}.
As special cases, we immediately get, in Section~\ref{sec4}, the Pokhozhaev conservation law for $I_2$ and the desired third-order conservation law for $I_3$, proving Theorem~\ref{T1}.
Finally, in Sections~\ref{sec6} and \ref{sec5}, we prove Lemmas~\ref{Teo3} and \ref{Teo2}
(and hence Theorem~\ref{Teo4}).


\hskip 0.3cm
\section{A quadratic form for a Liouville type equation} \label{LL}

Let us consider the Liouville-type (cf. \cite[\S6.2]{MU}) second-order ordinary differential equation
 \begin{equation}\label{lin2}
w_{tt} +\frac 1{q^2}\sp |\xi |^2\sp  w=0 \quad \text{for}  \quad  t\in [0,T),
\end{equation}
depending on the parameter $\xi \in \R^n\sp$,  with 
$$
q=q(t)\in C^4([0,T)), \quad q(t) \ne 0 \quad  \text{in}  \quad [0,T).
$$

Following \cite[\S\,2]{M1}, given a complex-valued solution 
$$
w= w(\xi, t)
$$ 
of \eqref{lin2}, we consider the 
quadratic form  $\spp {\pazocal E}= {\pazocal E}(\xi, t)$ defined by 
\begin{equation}
 \begin{aligned}\label{Q}
{\pazocal E} &:= \sum_{i=0}^1 \alpha_i|\xi|^{4-2i} \left(|w_t|^2  + \frac 1{q^2}|\xi|^2  |w|^2\right ) 
+\sum_{i=0}^1 \beta_i |\xi|^{4-2i}  \Re  (\overline{ w} \sp w_t)
 + \gamma_0  |\xi|^2 |w_t|^2,
 \end{aligned}
\end{equation}
where  $ \alpha_i=\alpha_i(t) , \beta_i=\beta_i(t)$  and $\gamma_0 =\gamma_0(t)$ 
are suitable  $C^1$ real  functions in $[0,T)$ to be determined later in order that
\begin{equation} \label{der27}
\frac{d}{dt} {\pazocal E}(\xi, t)=  \beta_1' |\xi|^2\Re (\overline{w} w_t).
\end{equation} 

To this aim we first compute
\beqs
\label{231}
\frac d{dt} \left [\alpha_0 |\xi|^4 \left (|w_t|^2  + \frac 1{q^2} |\xi|^2  |w|^2\right ) 
+ \beta_0 |\xi|^4  \Re  (\overline{ w}\sp w_t) + \gamma_0  |\xi|^2 |w_t|^2 \right ]
\eeqs
taking into account that $\spp w \spp $ solves \eqref{lin2} and hence
\beqs
\label{wt1}
\frac{d}{dt}|w_t|^2=w_{tt}\bar{w}_t+w_t \bar{w}_{tt}
=-\frac{1}{q^2}|\xi|^2w\bar{w}_t-\frac{1}{q^2}|\xi|^2w_t\bar{w}
=-\frac{2}{q^2}|\xi|^2\Re(\bar{w}w_t),
\eeqs
\beqs
\label{re1}
\frac{d}{dt}\Re(\bar w w_t)
=\Re(\bar{w}_tw_t+\bar{w}w_{tt})
=|w_t|^2-\frac{1}{q^2}|\xi|^2\Re(\bar{w}w)
=|w_t|^2-\frac{1}{q^2}|\xi|^2|w|^2.
\eeqs
Substituting  \eqref{wt1}, \eqref{re1} and $\frac{d}{dt}|w|^2=2\Re(\bar{w}w_t)$
into \eqref{231}, we get:
\begin{equation}
 \begin{aligned}
\frac d{dt} &\left [\alpha_0 |\xi|^4 \left (|w_t|^2  + \frac 1{q^2} |\xi|^2  |w|^2\right ) 
+ \beta_0 |\xi|^4  \Re  (\overline{ w}\sp w_t) + \gamma_0  |\xi|^2 |w_t|^2 \right ]
\\
&=\left [\left ( \frac{\alpha_0}{q^2} \right )'- \frac {\beta_0} {q^2}\right ]|\xi|^6 |w|^2 + \left [ \alpha_0' + \beta_0\right ]|\xi|^4 |w_t|^2 + 
\left [\beta_0' - 2 \frac{\gamma_0}{q^2}\right ] |\xi|^4 \Re (\overline{w}\sp w_t)
\\
&\quad + \gamma_0' |\xi|^2 |w_t|^2.
\end{aligned}
\end{equation}

Similarly,
\begin{equation}
  \begin{aligned}
\frac d{dt} &\left [\alpha_1 |\xi|^2 \left ( |w_t|^2  + \frac 1{q^2} |\xi|^2  |w|^2\right ) 
+ \beta_1 |\xi|^2  \Re  (\overline{ w}\sp w_t) \right ]
\\
&=\left [\left ( \frac{\alpha_1}{q^2} \right )'- \frac {\beta_1} {q^2}\right ]|\xi|^4 |w|^2 + \left [ \alpha_1' + \beta_1\right ]|\xi|^2 |w_t|^2 + 
\beta_1' |\xi|^2 \Re (\overline{w}\sp w_t).
\end{aligned}
\end{equation}

Let us now find
$\alpha_0, \beta_0, \gamma_0$ satisfying
\begin{equation} \label{S1}
{\begin{cases}
\displaystyle{\left ( \frac{\alpha_0}{q^2} \right )'- \frac {\beta_0} {q^2}=0}
 \\
 \displaystyle{\alpha_0' + \beta_0 =0 \phantom{\int^8_8} }
 \\
 \displaystyle{\beta_0' - 2 \frac{\gamma_0}{q^2} =0,}
 \end{cases}}
 \end{equation}
and then $\alpha_1, \beta_1$ such that
\begin{equation}\label {S2}
 {\begin{cases}
\displaystyle{\left ( \frac{\alpha_1}{q^2} \right )'- \frac {\beta_1} {q^2}=0}
 \\
 \displaystyle{\alpha_1' + \beta_1 = - \gamma_0'\sp,}
 \end{cases}}
 \end{equation}
so that the desired equality \eqref{der27} is satisfied.

From the first two equations of \eqref{S1} we find
\begin{equation}
\left ( \frac{\alpha_0}{q^2} \right )'- \frac {\beta_0} {q^2} = 
\frac 2 q  \left (\frac {\alpha_0} q \right )' =0.
\end{equation}
Hence, we have
\begin{equation} \label{c0}
\alpha_0 = C_0 \sp q, \quad \beta_0 = - C_0 \sp q' , \quad \gamma_0 = - \frac  {C_0} 2 \sp q^2 \sp q'',
\end{equation}
with $C_0\in \R$  arbitrary. Next, from system \eqref{S2}, we find
\begin{equation}
\left ( \frac{\alpha_1}{q^2} \right )'- \frac {\beta_1} {q^2}  = 
\frac 2 q  \left (\frac {\alpha_1} q \right )'  + \frac{\gamma_0'}{q^2}=0,
\end{equation}
and hence
\begin{equation}\label{der211}
\left (\frac {\alpha_1} q \right )' = - \frac{1} 2 \spi \frac{\gamma_0'}{q} 
=\frac{C_0}{4}(2q'q''+qq''')
\end{equation}
because of \eqref{c0}. Therefore
\begin{equation} \label{c2}
\alpha_1= \frac{C_0} 4 \sp \sp q  \left ( \frac {{q'}^2} 2 + q q''\right ) + C_1 \sp q, 
\end{equation}
with $C_1 \in \R$ arbitrary. Finally, from the second equation of \eqref{S2}, we have
\begin{equation} \label{c3}
 \begin{aligned} 
\beta_1 &= -\gamma_0' - \alpha_1' \\
&=\frac{C_0}{2}(q^2q'')'-\frac{C_0}{4}\left[q\left(\frac{{q'}^2}{2}+qq''\right)\right]'-C_1q'\\
&= -\frac{C_0}4 \left (  \frac {{q'}^3} 2 - q q'q'' - q^2 q'''\right ) -
C_1\sp q'.
\end{aligned}
\end{equation}
We have thus found $\alpha_0,\beta_0,\gamma_0,\alpha_1,\beta_1$ such that \eqref{der27}
is satisfied.


\hskip 0.2cm
\section{The quadratic form in the case of equation \eqref{KPR}} 
\label{sec3}

We now apply the results of the Section \ref{LL} to the Kirchhoff-Pokhozhaev equation \eqref{KPR}. 
Therefore, from now on, we assume that $\sp u\sp $ satisfies all the assumptions of Theorem~\ref{T1}.

First, since $\sp u\sp $ satisfies \eqref{reg2} and is a solution to
 \eqref{KPR}, from \eqref{reg} we easily get that
\begin{equation}
u \in C^k([0,T); H^{3-k}(\R^n)), \quad 0\le k \le 3 \sp ,
\end{equation}
and that
$$
\, t \,\mapsto \, \int_{\R^x} |\nabla u(x,t)|^2 dx \,
$$ 
is a $\sp C^4\sp $ function in $[0,T)$.

Taking the partial Fourier transform of $\sp u= u(x,t) \sp$ with respect to the  variable  $\sp x\sp $, i.e.
\begin{equation}
w(\xi,t) := (2\pi)^{-n/2} \int_{\R^n_x} u(x,t) e^{-i x \cdot \xi} \spp dx,  \quad \xi \in \R^n,
\end{equation}
and noting that, by Plancharel theorem,
\beqs
\label{32'}
\,\int_{\R^n} |\nabla u(x,t)|^2 dx= \int_{\R^n} |\xi|^2 |w(\xi,t)|^2 d\xi \spp ,
\eeqs
we have that $w$ is a complex valued solution of the equation
 \begin{equation}
 \label{KP3}
w_{tt} + \frac {|\xi|^2 \spp w} {\big ( a\int_{\R^n} |\xi |^2 |w |^2 d\xi+ b \big )^2 } 
\spp  = 0 ,\quad 
   t\in [0,T).
   \end{equation}
Next we consider, for $s=s(t)$ and 
$$
q= a \|\nabla u \|^2 + b = a\sp s +b
$$ 
as in \eqref{qs} and \eqref{reg2}, the quadratic form ${\pazocal E}(\xi, t)$ defined in \eqref{Q}. 
In the definition of $\spp {\pazocal E}(\xi, t)\spp $ we use the coefficients
$\alpha_0, \beta_0 $ and $\gamma_0$, defined as in \eqref{c0}, and $\alpha_1$, $\beta_1$ as in \eqref{c2}, \eqref{c3}. 

Finally, we introduce the functional 
 \begin{equation}
 E= E(t) := \int_{\R^n} {\pazocal E}(\xi, t)   d\xi \quad \text{for} \quad t \in [0,T).
 \end{equation}
From \eqref{qs} and \eqref{32'} we get
\begin{equation} \label{s1}
\int_{\R^n} |\xi|^2\Re (\overline{w} w_t) d\xi = \frac {1}{2}\spi s'\sp ,
\end{equation}
and hence formula \eqref{der27} for the derivative of ${\pazocal E}(\xi, t)$ gives
\begin{equation}\label{derE}
\frac{d}{dt} E=  \beta_1' \int_{\R^n} |\xi|^2\Re (\overline{w} w_t) d\xi = \frac {1}{2}\spi \beta_1'\spi s' \sp .
\end{equation}
Now, using the explicit definition of $\beta_1= \beta_1(t)$ given in \eqref{c3},  it is easy to show that 
the quantity $\beta_1' \spi s'$ is the  derivative of a polynomial in $\sp s, s', s''\sp $ and $\sp s'''\sp$. We have indeed:
\begin{equation}
\begin{aligned}
\beta_1' \spi s' &= \big ( \beta_1\spi s'\spi \big )' - \beta_1 \spi s'' 
\\
&= \big ( \beta_1\spi s'\spi \big )'  +\frac{C_0}4 \left (  \frac {{q'}^3} 2 - q q'q'' - q^2 q'''\right ) s'' +
C_1\sp q'\sp s''
\\
&= \big ( \beta_1\spi s'\spi \big )'  +\frac{C_0}4 \left (  \frac {a^3\, {s'}^3} 2 - a^2(as+b) s's'' - a(as+b)^2 s'''\right )\sp s'' +
C_1 \sp a\sp \sp s'\sp s''
\\
&= \big ( \beta_1\spi s'\spi \big )' 
 +\frac{C_0\sp a}8 \left (  \frac {a^2\sp {s'}^4} 4 - {(as+b)^2\spi {s''}^2} \right )' + 
\frac {C_1\sp a}{2}\spi \big ({s'}^2\big )'.
\end{aligned}
\end{equation}
Therefore, taking into account \eqref{derE}, we have proved that the quantity 
\begin{equation}
Q \spp := \spp  E- \frac {\beta_1\spi s'}{2}   -\frac{C_0\sp a }{16}  \left (  \frac {a^2\sp {s'}^4} 4- {(as+b)^2\spi {s''}^2} \right ) -
\frac {C_1 \sp a}{4}\spi {s'}^2
\end{equation}
remains constant in $[0,T)$.

Finally, let us compute  $\sp Q\sp$ explicitly, taking into account \eqref{s1} and the definitions 
of the coefficients $\alpha_0, \beta_0, \gamma_0$ and $\alpha_1, \beta_1$ given in \eqref{c0}, \eqref{c2} and \eqref{c3}. We have:
\begin{equation}
\begin{aligned}
\label{QQ}
Q \,=\,{}&
 C_0 \left ( q\spp {\| |\xi |^2 w_t\|^2} + \frac{\| |\xi |^3 w\|^2} q \right ) - C_0\spp q'
\int_{\R^n} |\xi|^4\Re (\overline{w} w_t) d\xi - \frac{C_0 }{2} q^2 q'' {\| |\xi | w_t\|^2}
\\
 &+ \left [\frac{C_0} 4 q\sp  \left ( \frac {{q'}^2} 2 + q q''\right ) + C_1q\right ]
  \left ( {\||\xi| w_t\|^2}  +  \frac{\||\xi|^2 w\|^2} {q^2}\right )
  \\
   &
  -\frac{C_0\sp a }{16}  \left (  \frac {a^2\sp {s'}^4} 4- {q^2\spi {s''}^2} \right ) -
\frac {C_1 \sp a}{4}\spi {s'}^2,
\end{aligned}
\end{equation}
where $\, C_0, C_1\in \R\,$ are arbitrary constants.


\hskip 0.2cm

\section{Conclusion of the Proof of Theorem\,\ref{T1}}
\label{sec4}

From \eqref{QQ}, we immediately deduce the conservation laws for the functional $I_2$ and $I_3$ as special cases of $Q$.

\hskip 0.1cm

\;\,i)\, If $\,C_0= 0$ and $\,C_1=1$, then
\begin{equation}
\label{41}
Q= 
  q\spp {\||\xi| w_t\|^2}  +  \frac{\||\xi|^2 w\|^2} {q} - \frac {a}{4}\spi {s'}^2= I_2,
\end{equation}
where the last equality is obtained by applying the Plancharel theorem to the
Pokhozhaev second order functional  \eqref{CL2}, written in the form 
\begin{equation}\label{I2bis}
I_2:= q {\|\nabla u_t\|^2} + \frac{\|\Delta u\|^2} q - \frac a 4 \spp {s'}^2,
\end{equation}
because of \eqref{s's''}.

\hskip 0.1cm

\;\, ii)\, If $\, C_0= 1$ and $\,C_1=0$, then  
\begin{equation}
\begin{aligned}\label{CL33}
Q&=  q\spp {\| |\xi |^2 w_t\|^2} + \frac{\| |\xi |^3 w\|^2} q  - \spp q'
\int_{\R^n} |\xi|^4\Re (\overline{w} w_t) d\xi - \frac{1 }{2} q^2 q'' {\| |\xi | w_t\|^2}
\\
 &\quad + \left [\frac{1} 4 q  \left ( \frac {{q'}^2} 2 + q q''\right ) \right ]
  \left ( {\||\xi| w_t\|^2}  +  \frac{\||\xi|^2 w\|^2} {q^2}\right )
  -\frac{ a }{16}  \left (  \frac {a^2\sp {s'}^4} 4- {q^2\spi {s''}^2} \right ) .
\end{aligned}
\end{equation}
Noting that $\spp q'' = a\sp s''\spp$ and
\begin{equation}
\label{44}
s'' = 2 \sp \||\xi| w_t\|^2 - \frac 2{q^2}\sp \||\xi|^2 w\|^2
\end{equation}
by  \eqref{s's''} (here and below we repeatedly apply Plancharel's theorem),
 we can write \eqref{CL33} as
\beqsn
\nonumber
Q=&&
q\spp \||\xi|^2w_t\|^2+\frac{\||\xi|^3w\|^2}{q}
-q'\int_{\R^n}|\xi|^4\Re(\bar{w}w_t)d\xi
-\frac12 q^2q''\||\xi|w_t\|^2\\
\nonumber
&&+\frac18 {q'}^2\left(q \spp \||\xi|w_t\|^2+\frac{\||\xi|^2w\|^2}{q}\right)
+\frac14q^2q''\||\xi|w_t\|^2
+\frac14q^2q''\frac{\||\xi|^2w\|^2}{q^2}\\
\nonumber
&&-\frac{a}{16}\frac{a^2{s'}^4}{4}
+\frac{a}{16}q^2\left(2 \sp \||\xi|w_t\|^2-\frac{2}{q^2}\||\xi|^2w\|^2\right)^2\\
\nonumber
=&&q\spp \||\xi|^2w_t\|^2+\frac{\||\xi|^3w\|^2}{q}
-q'\int_{\R^n}|\xi|^4\Re(\bar{w}w_t)d\xi
+\frac18 {q'}^2\left(q \spp \||\xi|w_t\|^2+\frac{\||\xi|^2w\|^2}{q}\right)
-\frac{a}{16}\frac{a^2{s'}^4}{4}\\
\nonumber
&&-\frac14aq^2 \sp \||\xi|w_t\|^2\left(2\||\xi|w_t\|^2-\frac{2}{q^2}\||\xi|^2w\|^2\right)
+\frac14aq^2\frac{\||\xi|^2w\|^2}{q^2}\left(2\||\xi|w_t\|^2-\frac{2}{q^2}\||\xi|^2w\|^2\right)\\
\nonumber
&&+\frac{1}{16}aq^2\left(2 \sp \||\xi|w_t\|^2-\frac{2}{q^2}\||\xi|^2w\|^2\right)^2
\\
\nonumber
=&&q \spp \||\xi|^2w_t\|^2+\frac{\||\xi|^3w\|^2}{q}
-q'\int_{\R^n}|\xi|^4\Re(\bar{w}w_t)d\xi
+\frac18 {q'}^2\left(q \spp \||\xi|w_t\|^2+\frac{\||\xi|^2w\|^2}{q}\right)
-\frac{a}{16}\frac{a^2{s'}^4}{4}\\
\nonumber
&&-\frac{1}{16}aq^2\left(2 \sp \||\xi|w_t\|^2-\frac{2}{q^2}\||\xi|^2w\|^2\right)
\left(4\sp \||\xi|w_t\|^2-4 \sp \frac{\||\xi|^2w\|^2}{q^2}-2 \sp \||\xi|w_t\|^2+\frac{2}{q^2}\||\xi|^2w\|^2\right).
\eeqsn
From \eqref{44} we finally get
 \begin{equation}
\begin{aligned}\label{CL34}
Q=&\,
   q\spp  {\| |\xi |^2 w_t\|^2} + \frac{\| |\xi |^3 w\|^2} q  -  q'
\int_{\R^n} |\xi|^4\Re (\overline{w} w_t) d\xi \\
&+ \frac{1} 8\spp  {q'}^2 \spp 
  \left ( { q \spp \||\xi| w_t\|^2}+  \frac{\||\xi|^2 w\|^2} {q}\right )
  -\frac{ a }{16}  \left (  \frac {a^2\sp {s'}^4} 4+ {q^2\spi {s''}^2} \right ),
\end{aligned}
\end{equation}
which is the third-order functional $I_3$ defined by the
 expression \eqref{I3}.


 \hskip 0.2cm
\section{Proof of Lemma\,\ref{Teo3}}
\label{sec6}

We note that, having $a\ne 0$ and $b > 0$, 
 the functional $I_1$
is (in any case) well defined  for $\spp \| \nabla u\|^2 < \frac b{|a|}\spp $ and we can rewrite the expression \eqref{I1KP} as
\begin{equation}\label{I11}
I_1
=  \|u_t\|^2 +  {\frac 1 b}\spp \frac {\| \nabla u\|^2}  {q},
\end{equation}
where $\sp q= a\|\nabla u\|^2 +b >0\sp $. 
Now, it is easy to see that 
\begin{equation}
I_1\leq \frac 1 {3|a|b}  \quad \Rightarrow \quad  \|\nabla u\|^2 \le \frac b {2 |a|}\sp .
\end{equation}
In fact, if $\spp I_1\leq \frac 1 {3|a|b}\,$, then \eqref{I11} gives
\beqsn
\|\nabla u\|^2\leq\frac{q}{3|a|}= \frac a{3|a|} \|\nabla u \|^2 + \frac b {3|a|}
\eeqsn
and hence
\beqsn
\frac23 \sp \|\nabla u\|^2\leq\left(1-\frac{a}{3|a|}\right)\|\nabla u\|^2\leq\frac{b}{3|a|}.
\eeqsn
Therefore $\sp \|\nabla u\|^2\leq \frac b {2|a|} \sp $ and hence we immediately get
\beqs
\label{stella}
\frac b2 \, \leq \,  q=a\|\nabla u\|^2+b \, \leq \, \frac32 \spp b \sp .
\eeqs
That is, the first statement in \eqref{condd} is satisfied. 

Finally, let us prove the second statement of \eqref{condd}. Taking into account the definition \eqref{lambd}
 of $\sp \lambda\sp$, we note that
\begin{equation}\label{st1}
|a|\lambda  \le |a| q I_1 + |a|b I_1 = |a| (q+b) I_1\sp .
\end{equation}
Therefore, if $\spp I_1 \le \frac 1 {3|a| b}\,$, from \eqref{stella} and
 \eqref{st1} we obtain that
\begin{equation}
|a| \lambda  \le 3 |a| b \spp I_1\sp .
\end{equation}
Hence, we find that
\begin{equation}
 I_1 \le \frac 1  {6|a| b}\;\; \Rightarrow \; \;  |a| \lambda \le  \frac 1 2 \sp .
\end{equation}

 
 \hskip 0.2cm
\section{Proof of Lemma\,\ref{Teo2}}
\label{sec5}

The proof of the estimate  \eqref{stima0} is almost immediate. From \eqref{s's''}, applying 
 Holder's inequality and taking into account that $q> 0$, we get
 \begin{equation}
 \label{stim1}
 \begin{aligned}
|s'| &\le 2 \int| \xi|^2 |w| |w_t| d\xi \\
&\leq2\left(\int_{\R^n}|\xi|^3|w|^2\right)^{1/2}\left(\int_{\R^n}|\xi||w_t|^2\right)^{1/2}\\
&\le 2\,  \| |\xi |  w \|^{\frac 1 2}\cdot
 \| |\xi |^2 w \|^{\frac 1 2} \cdot
 \|  w_t \|^{\frac 1 2}\cdot  \| |\xi | w_t \|^{\frac 1 2}\\
&
\le
\left ( q \spp \| |\xi | w_t \|^2 + \frac{  \| |\xi |^2 w \|^2 } q\right )^{ \frac 1 2}
 \left ( q \spp \|  w_t \|^2 + \frac{  \| |\xi | w \|^2 } q\right )^{ \frac 1 2}.
\end{aligned}
 \end{equation}
 Therefore, we deduce that
\begin{equation} \label{bound}
|a| \,\frac {{s'}^2} 4 \le \frac 1 2 \left ( q \spp \| |\xi | w_t \|^2 + \frac{  \| |\xi |^2 w \|^2 } q\right ),
\end{equation}
if $\, |a| \left ( q\|  w_t \|^2 + \frac{  \| |\xi | w \|^2 } q \right )  \le 2\spp$. 
The estimate \eqref{bound} applied to the expression \eqref{41} of $I_2$ gives
\beqsn
\frac 1 2 \left ( q\| |\xi | w_t \|^2 + \frac{  \| |\xi |^2 w \|^2 } q\right )
\leq I_2\leq \frac 32 \left ( q\| |\xi | w_t \|^2 + \frac{  \| |\xi |^2 w \|^2 } q\right )
\eeqsn
and hence \eqref{stima0}, since
\beqs
\label{lambda2}
\lambda=q \spp \|w_t\|^2+\frac{\||\xi|w\|^2}{q}
\eeqs
by \eqref{lambd}.
\hskip 0.2cm

Next let us come to the proof of \eqref{stima}. Again from \eqref{s's''} and Holder's inequality, we find:
 \begin{equation}
 \label{stim1'}
|s'| \le 2 \int| \xi|^2 |w| |w_t| d\xi 
\leq2\, \||\xi|^2w\|\cdot\|w_t\|
\le 2 \, \| |\xi |^3 w \|^{\frac 1 2} \cdot
 \| |\xi | w \|^{\frac 1 2}\cdot \|  w_t \|
 \end{equation}
 and also
  \begin{equation} \label{stim2}
|s'|\le  2  \| |\xi |^2 w_t \|^{\frac 1 2} \cdot
 \|  w_t \|^{\frac 1  2}\cdot \|| \xi |  w \|.
\end{equation}
Thus, having $q > 0$, from \eqref{stim1'} and \eqref{stim2} we deduce that
\begin{equation} \label{stim3}
\begin{aligned}
|s'| &\le  2 \, \| |\xi |^2 w_t \|^{\frac 1 4} \cdot \| |\xi |^3 w \|^{\frac 1 4} \cdot
 \|  w_t \|^{\frac 3  4}\cdot  \|| \xi |  w \|^{\frac 3 4}
 \\
 &\le  \left ( q \spp \| |\xi |^2 w_t \|^2 + \frac{  \| |\xi |^3 w \|^2 } q\right )^{ \frac 1 4}
 \left ( q \spp \|  w_t \|^2 + \frac{  \| |\xi | w \|^2 } q\right )^{ \frac 3 4}.
   \end{aligned}
\end{equation}
Moreover, again from Holder's inequality, we have
\begin{equation}
 \begin{aligned}
 \int| \xi|^4 |w| |w_t| d\xi 
 \leq \||\xi|^2w_t\|\cdot\||\xi|^2w\|
 \le   \| |\xi |^2 w_t \|\cdot   \| |\xi |^3 w \|^{\frac 1 2} \cdot
 \|| \xi |  w \|^{\frac 1 2}.
\end{aligned}
\end{equation}
Therefore, using \eqref{stim1'}, we obtain
\begin{equation}
\label{67}
 \begin{aligned}
|s'| \int| \xi|^4 |w| |w_t| d\xi & \le  2 \,  \| |\xi |^2 w_t \| \cdot \| |\xi |^3 w \| \cdot
 \|  w_t \|\cdot  \|| \xi |  w \|
 \\ 
 &\le {\frac  1 2 }
 \left ( q \spp \| |\xi |^2 w_t \|^2 + \frac{  \| |\xi |^3 w \|^2 } q\right )
 \left ( q \spp \|  w_t \|^2 + \frac{  \| |\xi | w \|^2 } q\right ).
\end{aligned}
\end{equation}
Similarly, by  \eqref{44} 
we find
 \begin{equation}
 \label{58}
 \begin{aligned}
q \spp |s''| &\le  2 \left (q\spp\||\xi| w_t\|^2 + \frac{ \| |\xi|^2 w \|^2} q \right )
\\
&\le  2 \left (q \spp \||\xi|^2 w_t\| \cdot \|w_t\| + \frac{ \| |\xi|^3 w \| \cdot \| |\xi| w \|} q \right )
\\
&\le  2 \left (q \spp \||\xi|^2 w_t\|^2  + \frac{ \| |\xi|^3 w \|^2 } q \right )^{\frac 1 2}
\left (q \spp \| w_t\|^2  + \frac{ \| |\xi| w \|^2 } q \right )^{\frac 1 2},
\end{aligned}
 \end{equation}
and, from $q'=as'$ and \eqref{stim1},
\begin{equation}
\label{69}
 \begin{aligned}
&{q'}^2 \left (q \spp \||\xi| w_t\|^2 + \frac{ \| |\xi|^2 w \|^2} q \right ) \le 
\\
&\leq a^2\left(q \spp \||\xi| w_t\|^2 + \frac{ \| |\xi|^2 w \|^2} q \right )^2
\left (q \spp \| w_t\|^2  + \frac{ \| |\xi| w \|^2 } q \right )\\
&\le  a^2 \left (q \spp \||\xi|^2 w_t\|^2  + \frac{ \| |\xi|^3 w \|^2 } q \right )
\left (q \spp \| w_t\|^2  + \frac{ \| |\xi| w \|^2 } q \right )^{2}.
\end{aligned}
 \end{equation}
 
In conclusion, taking into account the expression \eqref{lambda2} for $\lambda$,
from \eqref{67}, \eqref{69}, \eqref{stim3} and \eqref{58}, we can estimate the following expression which appear in the expression \eqref{CL34} of $Q$ (hence of $I_3$):
\beqsn
&&\left|-  q'
\int_{\R^n} |\xi|^4\Re (\overline{w} w_t) d\xi 
+ \frac{1} 8\spp  {q'}^2 \spp 
  \left ( { q \spp \||\xi| w_t\|^2}  +  \frac{\||\xi|^2 w\|^2} {q}\right )
  -\frac{ a }{16}  \left (  \frac {a^2\sp {s'}^4} 4+ {q^2\spi {s''}^2} \right )\right|\\
  \leq&&
\left ( {\frac{|a| \lambda} 2} +  {\frac{a^2 \lambda^2} 8} + {\frac{|a|^3 \lambda^3} {64}} + 
 {\frac{|a| \lambda} 4}\right ) 
 \left (q \, \||\xi|^2 w_t\|^2  + \frac{ \| |\xi|^3 w \|^2 } q \right ). 
\eeqsn
It is then clear, again from \eqref{CL34}, that the estimate  \eqref{stima} 
is valid if 
\begin{equation}\label{condl}
{\frac{|a| \lambda} 2} +  {\frac{a^2 \lambda^2} 8} + {\frac{|a|^3 \lambda^3} {64}} + 
 {\frac{|a| \lambda} 4} \le \frac 1 2 .
\end{equation}
It is then easy to see that \eqref{condl} is certainly true if  $\spp | a|\sp\lambda \le \frac  1 2\spp$.


\vspace*{4mm}
{\bf Acknowledgments.}
The first author is member of the Gruppo Nazionale per l'Analisi Ma\-te\-ma\-ti\-ca, la
Probabilit\`a e le loro Applicazioni (GNAMPA) of the Instituto Nazionale di Alta Ma\-te\-ma\-ti\-ca (INdAM) and was partially supported by the Italian Ministry of University and Research, under PRIN2022 (Scorrimento) ``Anomalies in partial differential equations and applications", code: 2022HCLAZ8\_002, CUP: J53C24002560006.


\end{document}